\documentclass[preprint,12pt]{elsarticle}

\newtheorem{theorem}{Theorem}[section]

\newtheorem{corollary}{Corollary}[section]
\newtheorem{proposition}{Proposition}[section]

\newcommand{\ignore}[1]{}{}

\usepackage{amssymb}
\usepackage{amsmath}
\numberwithin{equation}{section}
\allowdisplaybreaks

\usepackage{color}
\usepackage{soul}
\usepackage[dvipsnames]{xcolor}

\usepackage[colorlinks=true, linkcolor=blue, citecolor=blue]{hyperref}

\def\1{{{\mbox{${\rm{1\negthinspace\negthinspace I}}$}}}}

\renewcommand{\theequation}{\arabic{equation}}

\newcommand\beq{\begin{equation}}
\newcommand\eeq{\end{equation}}
\renewcommand{\theequation}{\thesection.\arabic{equation}}
\usepackage{ifthen}
\usepackage{xkeyval}
\usepackage{todonotes}
\begin{document}
\begin{frontmatter}
\title{Deviation inequalities for a supercritical branching process in a random environment}
\author{Huiyi  Xu}
\cortext[cor1]{\noindent Corresponding author. \\
\mbox{\ \ \ \ \ \ \ \ \ \ \  \ \ \ \ } . }
\address{Center for Applied Mathematics,
Tianjin University, Tianjin 300072, China. }

\begin{abstract}
\indent Let $\left \{ Z_{n}, n\ge 0  \right \}$ be a supercritical branching process in an independent and identically distributed random environment $\xi =\left ( \xi _{n}  \right )_{n\geq 0} $. In this paper, we get some  deviation inequalities for $\ln \left (Z_{n+n_{0} } / Z_{n_{0} } \right ).$ And some  applications are given  for constructing  confidence intervals.
\end{abstract}

\begin{keyword} Deviation inequalities; Branching processes; Random environment
\vspace{0.3cm}
\MSC primary 60J80; 60K37; secondary  60F10
\end{keyword}

\end{frontmatter}

\section{Introduction}

The branching process in a random environment (BPRE) is a generalization of the Galton-Watson process  by adding environment random variables.
It was first introduced by Smith and Wilkinson \cite{Smith}.
The  BPRE  can be described in the following form. Assume that $\xi =\left ( \xi _{0},\xi _{1},...   \right ) $ is a sequence of independently identically distributed (i.i.d.) random variables and $\xi _{n}$ stands for  the random environment at time $n$. Each  random variable $\xi _{n}$ corresponds to a probability law $p\left ( \xi _{n}  \right )=\left \{  p_{n}\left ( i  \right ):i\in \mathbb{N}  \right \} $ on $\mathbb{N}=\left \{ 0,1,... \right \}$, that is $\mathbb{P}\left ( \xi _{n} =i \right )= p_{n}\left ( i  \right ),  i\geq 0.$ Hence,  $  p_{n}\left ( i  \right ) $ is non-negative and satisfies $\sum_{i=0}^{\infty }p_{n}\left ( i  \right )=1 $.
In the random environment $\xi$, a branching process $\left \{ Z_{n}, n\ge 0  \right \} $ can be defined by the following equations:
\begin{center}
$Z_{0}  =  1$,\ \ \ \ \  $Z_{n+1}  =  \sum\limits_{i  =  1}^{Z_{n} } N_{n,i}$\ \ \mbox{for all}\ \ $n\ge 0$,
\end{center}
where $N_{n,i}$ represents the number of children of the $i$-th individual in generation $n$. Conditioned on the environment $\xi $,  the random variables $
\left \{N_{n,i}, n\ge 0, i\ge 1    \right \} $ are independent of each other  and  the random variables $\left \{N_{n,i},   i\ge 1    \right \}$ have a common law $p\left ( \xi _{n}  \right )$. In the sequel, denote by $\mathbb{P}_{\xi } $  the conditional probability when the environment $\xi$ is given, called the quenched law as usual. And $\mathbb{P}\left ( dx,d\xi  \right )=\mathbb{P_{\xi } }\left ( dx \right )\tau \left ( d\xi  \right ) $ stands for the total law of the process, called annealed law, where $\tau$ is the law of the environment $\xi$. The corresponding quenched and annealed expectations are represented by $\mathbb{E}_{\xi } $ and $\mathbb{E}$ respectively. For $n\ge 0 $, denote
\begin{center}
$m_{n}:=m_{n}\left ( \xi  \right )=\sum\limits_{i=0}^{\infty }ip_{n}\left ( i  \right )$ \ \ \ \mbox{and} \ \ \ $\Pi _{n}=\mathbb{E}_{\xi }Z_{n}=\prod\limits_{i=0}^{n-1}m_{i}$.
\end{center}
By the definition of expectation, it is easy to see  that $m_{n}=\mathbb{E}_{\xi }  N_{n,i}$ for each $i\ge 1 $.
The asymptotic behavior of $\log Z_n$ is crucially affected by the
associated random walk
\begin{eqnarray}
S_{n}=\ln \Pi _{n}=\sum\limits _{i=1}^{n} X_{i},\ \ \ n\ge 1. \nonumber
\end{eqnarray}
For simplicity, let
$$X=X_{1}=\ln m_{0},\ \ \ \mu =\mathbb{E}X\ \ \ \mbox{and} \ \ \ \sigma ^{2}=\mathbb{E}\left ( X-\mu  \right )^{2} .  $$
We call $\mu$ the criticality parameter.
According to the value   $ \mu  > 0, \mu  = 0$, or $\mu  < 0,$
the BPRE is respectively called supercritical, critical, or subcritical.

Because critical and subcritical BPRE's  will inevitably go extinct, the study of these two cases mainly focuses on the survival probability and conditional limit theorems for the branching processes, see, for instance,    Afanasyev et al. \cite{Afanasyev,Boinghoff}  and Vatutin \cite{Vatutin}. For the supercritical BPRE, a
number of researches have been focused on moderate and large deviations, see B{\"o}inghoff and Kersting \cite{Boinghoff2}, Bansaye and Berestycki \cite{Bansaye}, Huang and Liu \cite{Huang}, Kozlo \cite{Kozlov}, Nakashima \cite{Nakashima}, Bansaye and B{\"o}inghoff \cite{Bansaye2},  B{\"o}inghoff \cite{Boinghoff3} and  Wang and Liu \cite{WL17}.

In this paper, we assume that
\begin{center}
$p_{0}\left ( \xi _{0}  \right )=0$\ \ \ $\mathbb{P}$-a.s.\ \ \ \mbox{and}\ \ \ $0<\sigma ^{2}<\infty$,
\end{center}
which implies that the BPRE is supercritical, $Z_{n}\to \infty $ and the random walk $\left \{ S_{n} ,n\ge 0 \right \}$ is non-degenerate.
Under the  conditions:
$ \mathbb{E}  \frac{Z_1^{p}}{m_0}   < \infty$ for a constant $ p>1$ and
$ \mathbb{E} \exp\{t(X -\mu)\} < \infty $ for some $t$ in a neighborhood of $0$,  Grama et al. \cite{Grama}  have established the Cram\'{e}r moderate deviation expansion for the BPRE, which implies in particular that
for $0 \leq x = o(\sqrt{n} )$ as $n\rightarrow \infty$,
\begin{equation}\label{cramer}
\Bigg| \ln \frac{\mathbb{P}\big( \frac{ \ln  Z_{n}  - n \mu \ }{  \sigma \sqrt{n}} \geq x  \big)}{1-\Phi(x)} \Bigg|  \leq C    \frac{  1+x^3    }{  \sqrt{n}  },
\end{equation}
where  $C$ is  a positive  constant. See also Fan et al. \cite{Fan}   with more general conditions.
Asymptotic expansions, no matter how precise, do not diminish the need for probability inequalities valid for all $n, x$.
For the critical Galton-Watson process,  such type inequalities have been well studied by Nagaev \cite{NV}.
However, there are few papers on probability inequality for the BPRE.
 In order to fill this gap, we try to establish some deviation inequalities
 for the supercritical BPRE under  various moment conditions on $X$.

The paper is organized as follows. In Section 2, we present  our main results. In Section 3, some applications of the main results are discussed. The proofs of the main results are given in Section 4.
\section{Main results}
 To shorten notations, denote
$$Z_{n_{0},n }=\frac{\ln \frac{Z_{n+n_{0}} }{Z_{n_{0} } }-n\mu  }{\sigma \sqrt{n} } ,\ \ \ \ \ n_{0}, n \in \mathbb{N}.$$
In this section, we present  some deviation inequalities for $Z_{n_{0},n }$, under various moment conditions on $X$.

\subsection{Bernstein's inequality}
When $X$ satisfies Bernstein's condition, we have the following  Bernstein type inequality for $\ln \frac{Z_{n+n_{0}} }{Z_{n_{0} } }$. We refer to De la Pe\~{n}a \cite{D99} for similar results, where Bernstein type inequality for martingales is established.
\begin{theorem}\label{theorem 2.0}
Assume that there exists a positive  constant $H$ such that
\begin{eqnarray}\label{0gh}
\mathbb{E}(X -\mu)^{k}   \leq \frac12 k!H^{k-2} \mathbb{E} (X -\mu)^2  \    \ \ \textrm{for all}\   k\geq 2.
\end{eqnarray}
Then for all $ x > 0,$
\begin{equation}\label{fgh}
\mathbb{P}\left ( Z_{n_{0},n }\ge x  \right ) \leq 2\exp\bigg\{- \frac{  x^2  }{ 2(1+ 6\left ( 1+H \right ) \frac{x}{\sigma\sqrt{n}}  )  } \bigg\}.
\end{equation}
\end{theorem}

Condition (\ref{0gh}) is known as Bernstein's condition. It is known that Bernstein's condition is equivalent to Cram\'{e}r's condition: that is
$ \mathbb{E} \exp\{t(X -\mu)\} < \infty $ for some $t$ in a neighborhood of $0$, see Fan, Grama and Liu \cite{FGL12}.

From (\ref{fgh}), it is easy to see that for $0\leq x =o(\sqrt{n})$, the bound (\ref{fgh}) behaves as $2\exp\{-x^2/2\};$ while for $ x  \geq \sqrt{n}$,
it behaves as  $ 2\exp\{- c\, x \sqrt{n} \}$ for a constant $c>0$.

\subsection{Semi-exponential  inequality}
When $X$ has a  semi-exponential moment,   the following  theorem
holds. This  theorem can be compared to the corresponding results in Borovkov \cite{Bor} for partial
sums of independent random variables,   Dedecker et al. \cite{DDF19}  for Lipschitz functionals of composition of random functions, and Fan et al. \cite{X.Fan} for martingales.
\begin{theorem}\label{theorem 2.1}
	Assume
$  \mathbb{E}[ (X-\mu)^2 \exp\{ ((X-\mu)^+)^{\alpha} \}]< \infty  $
for some $\alpha \in \left ( 0,1 \right )$. Then  for all $ x > 0,$
\begin{eqnarray}
\mathbb{P}\left ( Z_{n_{0} ,n}\ge x  \right ) \le 3\exp\left \{ -\frac{x ^{2}}{8\left (  u+\big( \sigma \sqrt{n}  \right )^{-\alpha }x^{2-\alpha }  \big)   }  \right \},\label{0.1}
\end{eqnarray}
where
$$u= \frac{1}{ \sigma^2} \mathbb{E}[ (X-\mu)^2 \exp\{ ((X-\mu)^+)^{\alpha} \}]. $$
\end{theorem}

For moderate $0\leq x=o(n^{\alpha/(4-2\alpha)}) $, the
 bound (\ref{0.1}) is a sub-Gaussian bound and is of the order
\begin{eqnarray}
  3\exp\Big\{ -  \frac{x^2}{8u }   \Big\}.
\end{eqnarray}
For large  $x\geq n^{\alpha/(4-2\alpha)},$ the bound (\ref{0.1}) is a semi-exponential bound and is of the order
\begin{eqnarray}\label{constant}
 3  \exp\Big\{ -c\, x^\alpha  n^{\alpha/2}  \Big\}.
\end{eqnarray}
where $c$ does not depend on $x$ and $n$.
In particular,   inequality (\ref{0.1}) implies the following large deviation result:   there exists a positive constant $c$ such
that  for all $x >0$,
\begin{eqnarray}
  \mathbb{P}\left(\frac{1}{n} \ln \frac{Z_{n+n_{0}} }{Z_{n_{0} } }-\mu   \geq n x   \right)
   \leq 3  \exp\Big\{ -c\, x^\alpha  n^{\alpha}  \Big\},
\end{eqnarray}
where $c$ does not depend on $x$ and $n$.

\subsection{Fuk-Nagaev type bound}
When   $X  $ has an absolute moment  of order $p\geq 2 $, we have the following Fuk-Nagaev type inequality for  $Z_{n_{0} , n}.$
\begin{theorem}\label{theorem 2.2}
Let $p \geq 2$.  Assume that $\mathbb{E}  |X-\mu|^p    < \infty. $ Then  for  all $x>0,$
\begin{eqnarray}
\mathbb{P}\left ( Z_{n_{0} ,n}\ge x  \right ) \le  \exp\bigg\{- \frac{  x^2 }{ 2 V^2  } \bigg\} +  \frac{C_p}{n^{(p-2)/2}x^p },\label{4.1}
\end{eqnarray}
where
$$V^2= (p+2)^2 e^p\,  \ \ \ \ \  \textrm{and}\ \ \ \ \ \  C_p = 2^{p+1 }\Big( 1+ \frac{2}{p} \Big)^p \mathbb{E}  \Big| \frac{X-\mu}{\sigma} \Big|^p .$$
\end{theorem}

The last inequality  implies the following large deviation result:   there exists a positive constant $c$ such
that  for all $x >0$,
\begin{eqnarray}
  \mathbb{P}\left(\frac{1}{n} \ln \frac{Z_{n+n_{0}} }{Z_{n_{0} } }-\mu   \geq n x   \right)
   \leq   \exp\bigg\{-\, \frac{x^2}{c} n \bigg\} +  \frac{c }{n^{p-1}x^p }.
\end{eqnarray}
 Thus for any $x>0$,
$$
 \mathbb{P}\left(\frac{1}{n} \ln \frac{Z_{n+n_{0}} }{Z_{n_{0} } }-\mu   \geq n x   \right)=O\Big(\frac{1}{ n^{p-1}} \Big) \,
$$
as $n\rightarrow \infty.$ Note that the last  equality is optimal under the stated condition.

\subsection{von Bahr-Esseen type bound}
When the  random variable $X $ has an absolute moment  of order $p\in \left ( 1,2 \right ] $, we have the following von Bahr-Esseen inequality.  Notice that in the next theorem, the variance of $X$ may not exist. Thus we consider the large deviation inequality for  $\ln \frac{Z_{n+n_{0} } }{Z_{n_{0} } }$  instead of $Z_{n_{0} ,n}$.
\begin{theorem}\label{theorem 2.5}
	Let $p\in \left ( 1,2 \right ] $. Assume that $\mathbb{E}  |X-\mu |^p    < \infty.$  Then for all $x>0$,
\begin{eqnarray}\label{sdfs}
\mathbb{P} \left (\frac{1}{n} \ln \frac{Z_{n+n_{0}} }{Z_{n_{0} } }-\mu  \ge x\right )\le\frac{C_{p} }{x^{p}n^{p-1}} ,
\end{eqnarray}
where
$$ C_p =  2^{p+1} \mathbb{E}|X-\mu  |^p+ (2p)^pe^{-p}   .$$
\end{theorem}

From the inequality (\ref{sdfs}), it is easy to see that for any $x>0,$
\begin{eqnarray}
\mathbb{P} \left (\frac{1}{n} \ln \frac{Z_{n+n_{0}} }{Z_{n_{0} } }-\mu  \ge x\right ) =O\Big(\frac{1}{n^{p-1}} \Big),\ \ \ \ \ n\rightarrow \infty.
\end{eqnarray}
The last convergence rate is the best possible under the stated condition,   see von Bahr and Esseen \cite{von Bahr}.

\subsection{Hoeffding type bound}
When    $X $  is bounded  from above, we obtain the following Hoeffding type inequality for  $Z_{n_{0},n }.$ We refer to Fan, Grama and Liu \cite{X.Fan1} for similar results, where Bernstein type inequality for martingales is established.
\begin{theorem}\label{theorem 2.6}
Assume that there exists a positive  constant $H$ such that
$$ X \leq \mu+ H.$$
Then for all $0<x\leq\frac{\sigma \sqrt{n} }{2} $,
\begin{eqnarray}
\mathbb{P}\left ( Z_{n_{0},n }\ge x  \right ) &\le & 2 \exp \left \{ -\frac{x}{2H } \left [ \left ( 1+\frac{2\sigma \sqrt{n}}{Hx }  \right )\ln \left ( 1+\frac{Hx }{2\sigma \sqrt{n}}  \right )  -1 \right ] \right \}\nonumber \\
&\le & 2 \exp \left \{ -\frac{x^{2} }{8\left ( 1+\frac{Hx}{6\sigma \sqrt{n} }  \right ) }  \right \};\label{10}
\end{eqnarray}
and for $x>\frac{\sigma \sqrt{n} }{2}  $,
\begin{eqnarray}
\mathbb{P}\left ( Z_{n_{0},n }\ge x  \right )& \le& \exp \left \{ -\frac{x}{2H } \left [ \left ( 1+\frac{2\sigma \sqrt{n}}{Hx }  \right )\ln \left ( 1+\frac{Hx }{2\sigma \sqrt{n}}  \right )  -1 \right ] \right \}+\exp\left \{- \frac{x\sigma \sqrt{n} }{2}  \right \}\nonumber\\
&\le& \exp \left \{ -\frac{x^{2} }{8\left ( 1+\frac{Hx}{6\sigma \sqrt{n} }  \right ) }  \right \}+\exp\left \{- \frac{x\sigma \sqrt{n} }{2}  \right \}.\label{11}
\end{eqnarray}
\end{theorem}

\subsection{Rio type bound}
When   $X  $ is bounded, we obtain the following Rio type inequality for  $\ln \frac{Z_{n+n_{0}} }{Z_{n_{0} } }.$
\begin{theorem}\label{theorem 2.8}
Assume that there exists two positive constants $H_{1}$ and $H_{2}$ such that
$$ H_{1}  \le X-\mu  \le H_{2} . \ \ \ \   $$
Then   for all $x\in [0,2\left ( H_{2}-H_{1}   \right )) $,
\begin{eqnarray}
\mathbb{P} \left (\frac{1}{n} \ln \frac{Z_{n+n_{0}} }{Z_{n_{0} } }-\mu  \ge x\right )&\le&\exp \Big \{ -n\max \left ( \psi _{1}\left ( x \right ), \psi _{2}\left ( x \right )   \right ) \Big \} +\exp \left \{ -\frac{1}{2}nx  \right \}\nonumber \\
&\le&\left ( 1-\frac{x}{2\left ( H_{2}-H_{1} \right )}  \right )^{\frac{nx}{H_{2}-H_{1}} \left ( 1-\frac{x }{4\left ( H_{2}-H_{1} \right ) }  \right ) }   +\exp \left \{ -\frac{1}{2}nx  \right \},\nonumber\ \ \ \
\end{eqnarray}
where $$\psi _{1}\left ( x \right )=\frac{x^{2} }{2\left ( H_{2}-H_{1}   \right )^{2}  } +\frac{x^{4} }{36\left ( H_{2}-H_{1}   \right )^{4}},\ \ \ \  x> 0,$$ and
\begin{eqnarray}
  \psi _{2}\left ( x \right )=\left ( \frac{x^{2} }{4\left ( H_{2}-H_{1}   \right )^{2}  }-\frac{x}{H_{2}-H_{1}}  \right )\ln \left ( 1-\frac{x}{2\left ( H_{2}-H_{1}   \right )}  \right ),\  \ x\in [0,2\left ( H_{2}-H_{1}   \right )). \label{0.8}
\end{eqnarray}
\end{theorem}

By Rio's remark \cite{Rio}, for all $x$ in $\left [ 0,1 \right ] $, we have
 $$\left ( 1-x \right )^{nx\left ( 2-x \right ) }\le \exp\left \{ -2nx^{2}  \right \} ,$$
 which leads to, for all $x\in [0,2\left ( H_{2}-H_{1}   \right )) $,
 $$\left ( 1-\frac{x}{2\left ( H_{2}-H_{1} \right )}  \right )^{\frac{nx}{H_{2}-H_{1}} \left ( 1-\frac{x }{4\left ( H_{2}-H_{1} \right ) }  \right ) }\le\exp \left \{ -\frac{nx^{2} }{2\left ( H_{2}-H_{1}   \right )^{2}   }  \right \}. $$
So we get the following corollary, a  simple consequence of the Rio type inequality.
\begin{corollary}\label{corollary 2.7}
Assume the condition of Theorem \ref{theorem 2.8}.
Then   for all $x\ge 0 $,
\begin{eqnarray}
	\mathbb{P} \left (\frac{1}{n} \ln \frac{Z_{n+n_{0}} }{Z_{n_{0} } }-\mu  \ge x\right )\le\exp \left \{ -\frac{nx^{2} }{2\left ( H_{2}-H_{1}   \right )^{2}   }  \right \}+\exp \left \{ -\frac{1}{2}nx  \right \} .\label{2.0}
\end{eqnarray}
\end{corollary}

When $0\le x\le \left ( H_{2}-H_{1}   \right )^{2} $, the second term in the right hand side of (\ref{2.0}) is less than the first one. Thus
we have the following  Azuma-Hoeffding  inequality: for all $0\le x\le \left ( H_{2}-H_{1}   \right )^{2} ,$
\begin{eqnarray}\label{1.0}
\mathbb{P} \left (\frac{1}{n} \ln \frac{Z_{n+n_{0}} }{Z_{n_{0} } }-\mu  \ge x\right )\le 2\exp \left \{ -\frac{nx^{2} }{2\left ( H_{2}-H_{1}   \right )^{2}   }  \right \}.
\end{eqnarray}

\section{Application to construction of confidence intervals}
Deviation inequalities can be applied to establishing confidence intervals for the criticality parameter $\mu$ in terms of $Z_{n_{0} } $, $Z_{n+n_{0}} $ and $n$, or to preview $Z_{n+n_{0}} $ in terms of $Z_{n_{0} } $, $\mu$ and $n$.
\subsection{Construction of confidence intervals for $\mu$}
When $Z_{n_{0} } $, $Z_{n+n_{0}} $ and $\sigma ^{2} $ are known, we can use Theorem \ref{theorem 2.0} to estimate $\mu$.
\begin{proposition}\label{proposition 3.1}
Assume that there exists a positive  constant $H$ such that
$$\mathbb{E}(X -\mu)^{k}   \leq \frac12 k!H^{k-2} \mathbb{E} (X -\mu)^2  \    \ \ \textrm{for all}\   k\geq 3.$$
Let $\delta _{n}\in \left ( 0,1 \right ) $ and
$$\Delta _{n}=\frac{6\left ( 1+H \right ) }{n} \ln \left ( 2/\delta _{n}  \right )+\sqrt{\frac{36\left ( 1+H \right ) ^{2} }{n^{2} } \ln ^{2} \left ( 2/\delta _{n}  \right )+\frac{2}{n}\sigma ^{2}\ln \left ( 2/\delta _{n}  \right )   }.$$
Then $\left [ A_{n},+\infty    \right ) $, with
$$ A_{n}=\frac{1}{n}\ln \left ( \frac{Z_{n_{0}+n } }{Z_{n_{0} } }  \right )-\Delta_{n}, $$
is a $1-\delta _{n}$ confidence interval for $\mu $.
\end{proposition}
\textit{Proof.} By Theorem \ref{theorem 2.0}, we have
\begin{eqnarray}\label{9.0}
\mathbb{P}\left ( Z_{n_{0},n }> \frac{x\sqrt{n} }{\sigma }   \right )\le 2\exp \left \{ -\frac{nx^{2} }{2\left ( \sigma ^{2}+6\left ( 1+H \right ) x   \right ) }  \right \} .
\end{eqnarray}
Let $\Delta _{n}$ be the positive solution of the following equation
\begin{eqnarray}\label{9.5}
2\exp \left \{ -\frac{nx^{2} }{2\left ( \sigma ^{2}+6\left ( 1+H \right ) x   \right ) }  \right \}=\delta _{n}.
\end{eqnarray}
Then
\begin{eqnarray}\label{9.6}
\Delta _{n}=\frac{6\left ( 1+H \right ) }{n} \ln \left ( 2/\delta _{n}  \right )+\sqrt{\frac{36\left ( 1+H \right ) ^{2} }{n^{2} } \ln ^{2} \left ( 2/\delta _{n}  \right )+\frac{2}{n}\sigma ^{2}\ln \left ( 2/\delta _{n}  \right )   }.
\end{eqnarray}
By \eqref{9.0}, we obtain
\begin{eqnarray}\label{9.3}
\mathbb{P}\left ( Z_{n_{0},n }> \frac{\Delta _{n}\sqrt{n} }{\sigma }   \right )=\mathbb{P}\left ( \mu < \frac{1}{n}\ln \left ( \frac{Z_{n_{0}+n } }{Z_{n_{0} } }  \right ) -\Delta _{n}  \right ) \le \delta _{n}.\nonumber
\end{eqnarray}
The last inequality implies that
\begin{eqnarray}
\mathbb{P}\left ( \mu \ge \frac{1}{n}\ln \left ( \frac{Z_{n_{0}+n } }{Z_{n_{0} } }  \right ) -\Delta _{n}  \right )\ge 1- \delta _{n}.\nonumber
\end{eqnarray}
This completes the proof of Proposition \ref{proposition 3.1}.

Similarly, when $X$ is bounded, we have the following estimate for $\mu$.
\begin{proposition}\label{proposition 3.2}
Assume that there exists positive constants $H_{1}$ and $H_{2}$ such that
$$ H_{1}  \le X-\mu  \le H_{2} .$$
Let $\delta _{n}\in \left [2\exp \left \{ -\frac{n}{2} \left ( H_{2}-H_{1}  \right )^{2}   \right \},1   \right ]  $ and
$$\Delta _{n}=\left (  H_{2}-H_{1}\right ) \sqrt{\frac{2 }{n}\ln \left ( 2/\delta _{n}  \right )  } .$$
Then $\left [ A_{n},+\infty    \right ) $, with
$$ A_{n}=\frac{1}{n}\ln \left ( \frac{Z_{n_{0}+n } }{Z_{n_{0} } }  \right )-\Delta_{n}, $$
is a $1-\delta _{n}$ confidence interval for $\mu $.
\end{proposition}
\textit{Proof.} By Corollary \ref{corollary 2.7}, we have
\begin{eqnarray}\label{9.1}
\mathbb{P} \left (\frac{1}{n} \ln \frac{Z_{n+n_{0}} }{Z_{n_{0} } }-\mu  > x\right )\le 2\exp \left \{ -\frac{x^{2} }{2n^{-1}\left ( H_{2}-H_{1}   \right )^{2}   }  \right \}.
\end{eqnarray}
Let $\Delta _{n}$ be the positive solution of the following equation
\begin{eqnarray}\label{9.7}
\delta _{n} = 2\exp \left \{ -\frac{nx^{2} }{2\left ( H_{2}-H_{1}   \right )^{2} }  \right \}.
\end{eqnarray}
Then
\begin{eqnarray}\label{9.8}
\Delta _{n} =\left ( H_{2}-H_{1} \right )  \sqrt{\frac{2}{n}\ln \left ( 2/\delta _{n}  \right )  }   .
\end{eqnarray}
It is easy to see that
\begin{eqnarray}
 \mathbb{P}\left ( \mu \ge \frac{1}{n}\ln \left ( \frac{Z_{n_{0}+n } }{Z_{n_{0} } }  \right ) -\Delta _{n}   \right )\ge 1- \delta _{n}.\nonumber
\end{eqnarray}
This completes the proof of Proposition \ref{proposition 3.2}.
\subsection{Construction of confidence intervals for $Z_{n+n_{0}} $}
When the parameter $\mu$ and $Z_{n_{0} } $ are known, we can use Theorem \ref{theorem 2.0} to preview $Z_{n+n_{0}} $.
\begin{proposition}\label{proposition 3.3}
Assume that there exists a positive  constant $H$ such that
$$\mathbb{E}(X -\mu)^{k}   \leq \frac12 k!H^{k-2} \mathbb{E} (X -\mu)^2  \    \ \ \textrm{for all}\   k\geq 3.$$
Let $\delta _{n}\in \left ( 0,1 \right ) $  and
$$\Delta _{n}=\frac{6\left ( 1+H \right ) }{n} \ln \left ( 2/\delta _{n}  \right )+\sqrt{\frac{36\left ( 1+H \right ) ^{2} }{n^{2} } \ln ^{2} \left ( 2/\delta _{n}  \right )+\frac{2}{n}\sigma ^{2}\ln \left ( 2/\delta _{n}  \right )   }.$$
Then $\left [ 1,A_{n}  \right ]$, with
$$A_{n}=Z_{n_{0} } \exp \left \{ n\left ( \mu +\Delta_{n} \right )  \right \} ,$$
is a $1-\delta _{n}$ confidence interval for $Z_{n+n_{0}} $.
\end{proposition}
\textit{Proof.} With arguments similar to that of \eqref{9.0}-\eqref{9.6}, we have
\begin{eqnarray}\label{0.3}
\mathbb{P}\left ( Z_{n_{0},n }> \frac{\Delta _{n}\sqrt{n} }{\sigma }   \right )=\mathbb{P}\left (Z_{n_{0}+n } >  Z_{n_{0} }\exp \left \{ n\left ( \mu +\Delta _{n} \right )  \right \}   \right ) \le \delta _{n}.\nonumber
\end{eqnarray}
Then
\begin{eqnarray}
\mathbb{P}\left (Z_{n_{0}+n }\le  Z_{n_{0} }\exp \left \{ n\left ( \mu +\Delta _{n}  \right )  \right \}   \right )\ge 1- \delta _{n}.\nonumber
\end{eqnarray}
This completes the proof of Proposition \ref{proposition 3.3}.

Similarly, when $X$ is bounded , the parameter $\mu$ and $Z_{n_{0} } $ are known, we can use Corollary \ref{corollary 2.7} to preview $Z_{n+n_{0}} $.
\begin{proposition}\label{proposition 3.4}
Assume that there exists positive constants $H_{1}$ and $H_{2}$ such that
$$ H_{1}  \le X-\mu  \le H_{2} .$$
Let $\delta _{n}\in \left [2\exp \left \{ -\frac{n}{2} \left ( H_{2}-H_{1}  \right )^{2}   \right \},1   \right ]  $ and
$$\Delta _{n}=\left (  H_{2}-H_{1}\right ) \sqrt{\frac{2 }{n}\ln \left ( 2/\delta _{n}  \right )  } .$$
Then $\left [ 1,A_{n} \right ] $, with
$$A_{n}=Z_{n_{0} } \exp \left \{ n\left ( \mu +\Delta_{n} \right )  \right \} ,$$
is a $1-\delta _{n}$ confidence interval for $Z_{n+n_{0}} $.
\end{proposition}
\textit{Proof.} Again by arguments similar to that of \eqref{9.1}-\eqref{9.8}, we get
\begin{eqnarray}
\mathbb{P} \left (\frac{1}{n} \ln \frac{Z_{n+n_{0}} }{Z_{n_{0} } }-\mu  > \Delta _{n}\right )= \mathbb{P}\left (Z_{n_{0}+n }>  Z_{n_{0} }\exp \left \{ n\left ( \mu +\Delta _{n}  \right )  \right \}   \right )\le \delta _{n}. \nonumber
\end{eqnarray}
It is easy to see that
\begin{eqnarray}
\mathbb{P}\left (Z_{n_{0}+n }\le  Z_{n_{0} }\exp \left \{ n\left ( \mu +\Delta _{n}  \right )  \right \}   \right )\ge 1- \delta _{n}.\nonumber
\end{eqnarray}
This completes the proof of Proposition \ref{proposition 3.4}.
\section{Proofs of Theorems }
Denote by
\begin{eqnarray}
W_{n}=\frac{Z_{n} }{\Pi _{n} } , \ \ \ n\ge 0, \label{1.1}
\end{eqnarray}
the normalized population size. As we all know, the sequence $\left ( W_{n}  \right )_{n\ge 0}$ is a positive martingale both under the quenched law $\mathbb{P}_{\xi }$ and under the annealed law $\mathbb{P} $ with respect to the natural filtration
$$\mathcal{F}_{0}=\sigma \left \{ \xi  \right \},\ \mathcal{F}_{n}=\sigma \left \{ \xi ,N_{k,i},0\le k\le n-1,i\ge 1  \right \} ,\ \ \ n\ge 1.$$
According to Doob's  convergence theorem and Fatou's lemma, the limit $W=\lim_{n\rightarrow\infty} W_{n} $ exists $\mathbb{P}$-a.s. and $\mathbb{E}W\le 1 $. Evidently, formula \eqref{1.1} implies the following decomposition:
\begin{eqnarray}
\ln Z_{n}=\sum\limits_{i=1}^{n}X_{i}+\ln W_{n}, \label{1.2}
\end{eqnarray}
where $X_{i}=\ln m_{i-1} \left ( i\ge 1 \right )$ are i.i.d. random variables depending only on the environment $\xi $. Consequently, the asymptotic behavior of $\ln Z_{n} $ is primarily affected by the associated random walk $S_n=\sum\limits _{i=1}^{n} X_{i}$.
In the sequel, we denote
\begin{eqnarray}
\eta _{n,i}=\frac{X_{i}-\mu}{\sigma \sqrt{n} }  ,\ \ \ i=1,...,n_{0}+n, \nonumber
\end{eqnarray}
and
$$W_{n_{0},n }=\frac{W_{n_{0}+n } }{W_{n_{0} } } .$$
Then it is easy to see that $\sum_{i=1}^{n}\mathbb{E}\eta_{n,n_{0}+i }^{2}=1 $ and $\sum\limits _{i=1}^{n} \eta _{n,n_{0}+i }$ is a sum of i.i.d. random variables.
\subsection{Proof of Theorem \ref{theorem 2.0}}
We first give a proof of the inequality for $0\le x<\frac{\sigma \sqrt{n} }{2}. $ Clearly, it holds for all $x\ge 0 ,$
\begin{eqnarray}
\mathbb{P}\left ( Z_{n_{0},n }\ge x  \right )& = & \mathbb{P}\left ( \sum\limits_{i=1}^{n}\eta _{n,n_{0}+i } +\frac{\ln W _{n_{0},n }}{\sigma \sqrt{n} } \ge x   \right )  \nonumber\\
& \le & I_{1}+I_{2},\label{3.1}
\end{eqnarray}
where\\
\begin{eqnarray}
I_{1}=\mathbb{P} \left ( \sum\limits_{i=1}^{n}\eta _{n,n_{0}+i } \ge \left ( x-\frac{x^{2} }{\sigma \sqrt{n} }  \right )  \right )  \  \ \mbox{and}\  \ I_{2}= \mathbb{P} \left ( \frac{\ln W _{n_{0},n }   }{\sigma \sqrt{n} }\ge \frac{x^{2} }{\sigma \sqrt{n} }  \right ).\label{3.2}
\end{eqnarray}
Next, we give some estimations for $I_{1}$ and $I_{2}$. Using Bernstein's inequality \cite{Bernstein} for i.i.d. random variables, we obtain for all $0\le x<\frac{\sigma \sqrt{n} }{2},$
\begin{eqnarray}
I_1 &\leq&   \exp\bigg\{- \frac{  x^2(1- \frac{  x}{\sigma \sqrt{n}} )^2 }{ 2(1+ \frac{H}{\sigma\sqrt{n}} x(1- \frac{x}{\sigma \sqrt{n}} ))  } \bigg\} \nonumber  \\
&\leq&  \exp\bigg\{- \frac{  x^2  }{ 2(1+ 6\left ( 1+H \right ) \frac{x}{\sigma\sqrt{n}}  )  } \bigg\} .\label{2.3}
\end{eqnarray}
By Markov's inequality and the fact that $\mathbb{E}W_{n}=1 $, we have for all $0\le x<\frac{\sigma \sqrt{n} }{2},$
\begin{eqnarray}
I_{2}&=&\mathbb{P}\left ( W_{n_{0},n } \ge \exp\left \{ x^{2}  \right \}  \right ) \nonumber \\
&\le& \exp\left \{ -x^{2}  \right \}\mathbb{E}W_{n_{0},n }=\exp\left \{ -x^{2}  \right \}\nonumber \\
&\le& \exp\bigg\{- \frac{  x^2  }{ 2(1+ 6\left ( 1+H \right ) \frac{x}{\sigma\sqrt{n}}  )  } \bigg\}.\label{2.4}
\end{eqnarray}
Combining \eqref{3.1}, \eqref{2.3} and \eqref{2.4}, we obtain for all $0\le x<\frac{\sigma \sqrt{n} }{2}$,
$$\mathbb{P}\left ( Z_{n_{0},n }\ge x  \right )\le  2\exp\bigg\{- \frac{  x^2  }{ 2(1+ 6\left ( 1+H \right ) \frac{x}{\sigma\sqrt{n}}  )  } \bigg\}.$$
When $x>\frac{\sigma \sqrt{n} }{2}$, it holds
\begin{eqnarray}
\mathbb{P}\left ( Z_{n_{0},n }\ge x  \right )\le I_{3}+I_{4} ,\label{3.5}
\end{eqnarray}
where
\begin{eqnarray}
I_{3}=\mathbb{P}\left ( \sum\limits_{i=1}^{n}\eta _{n,n_{0}+i }\ge \frac{x}{2}    \right ) \ \ \mbox{and} \ \ I_{4}=\mathbb{P}\left ( \frac{\ln W_{n_{0},n }    }{\sigma \sqrt{n} }\ge \frac{x}{2}   \right ).\label{3.6}
\end{eqnarray}
Again by Bernstein's inequality for i.i.d.\ random variables,   we obtain for all $x>\frac{\sigma \sqrt{n} }{2}$,
 \begin{eqnarray}
I_3 &\leq&   \exp\bigg\{- \frac{  (x/2)^2  }{ 2(1+ \frac{H}{\sigma\sqrt{n}} \frac x 2 )  } \bigg\}
 \leq   \exp\bigg\{- \frac{  x^2  }{ 2(1+ 6\left ( 1+H \right ) \frac{x}{\sigma\sqrt{n}}  )  } \bigg\} .
\end{eqnarray}
Again by  Markov's inequality and the fact that $\mathbb{E}W_n=1$,  we have for all $x>\frac{\sigma \sqrt{n} }{2}$,
\begin{eqnarray}
I_{4}&=&\mathbb{P}\left ( W_{n_{0},n }\ge \exp\left \{ \frac{x\sigma \sqrt{n} }{2}  \right \}   \right )\le \exp\left \{- \frac{x\sigma \sqrt{n} }{2}  \right \}\mathbb{E}W_{n_{0},n }\nonumber \\
&=&\exp\left \{- \frac{x\sigma \sqrt{n} }{2}  \right \}\nonumber \\
&\le &\exp\bigg\{- \frac{  x^2  }{ 2(1+ 6\left ( 1+H \right ) \frac{x}{\sigma\sqrt{n}}  )  } \bigg\}. \label{2.8}
\end{eqnarray}
Thus, for all $x>\frac{\sigma \sqrt{n} }{2}$,
$$\mathbb{P}\left ( Z_{n_{0},n }\ge x  \right )  \leq 2 \exp\bigg\{- \frac{  x^2  }{ 2(1+ 6\left ( 1+H \right ) \frac{x}{\sigma\sqrt{n}}  )  } \bigg\}. $$
This completes the proof of Theorem \ref{theorem 2.0}.  \hfill\qed
\subsection{Proof of Theorem \ref{theorem 2.1}}
 From \eqref{3.1} and \eqref{3.2}, using the inequality of Fan, Grama and Liu \cite{X.Fan} for i.i.d.\ random variables, we get for all $0\le x<\frac{\sigma \sqrt{n} }{2},$
\begin{eqnarray}
I_{1}& =& \mathbb{P}\left ( \sum\limits _{i=1}^{n} \left ( X _{n_{0}+i }-\mu  \right )\ge \sigma \sqrt{n}\left ( x-\frac{x^{2} }{\sigma \sqrt{n} }  \right )    \right )\nonumber \\
& \le  & 2\exp\left \{ -\frac{\left ( \sigma \sqrt{n}\left ( x-\frac{x^{2} }{\sigma \sqrt{n} }  \right )   \right )^{2}  }{2\left ( u_{n} + \left ( \sigma \sqrt{n}\left ( x-\frac{x^{2} }{\sigma \sqrt{n} }  \right )   \right )^{2-\alpha }  \right ) }  \right \}\nonumber \\
& = & 2\exp\left \{ -\frac{x ^{2}\left ( 1- \frac{x }{\sigma \sqrt{n} }\right )^{2}   }{2\left ( \frac{u_{n}}{\sigma ^{2}n } +\left ( \sigma \sqrt{n}  \right )^{-\alpha } x^{2-\alpha }  \left ( 1-\frac{x}{\sigma \sqrt{n} }  \right )^{2-\alpha }      \right ) }  \right \}\nonumber \\
& \le  & 2\exp\left \{ -\frac{x ^{2}}{8\left (  u+\left ( \sigma \sqrt{n}  \right )^{-\alpha }x^{2-\alpha }  \right )   }  \right \},\label{3.3}
\end{eqnarray}
where $u_{n} =n \mathbb{E}[ (X-\mu)^2 \exp\{ ((X-\mu)^+)^{\alpha} \}]$ and $I_{i}(i=1,2,3,4)$ are as defined in the proof of Theorem \ref{theorem 2.0}. By Markov's inequality and the fact that $\mathbb{E}W_{n}=1 $, we have for all $0\le x<\frac{\sigma \sqrt{n} }{2},$
\begin{eqnarray}
I_{2}\le \exp\left \{ -x^{2}  \right \}.\label{3.4}
\end{eqnarray}
Combining \eqref{3.1}, \eqref{3.3} and \eqref{3.4}, we obtain for all $0\le x<\frac{\sigma \sqrt{n} }{2},$
\begin{eqnarray}
\mathbb{P}\left ( Z_{n_{0},n } \ge x \right )& \le & 2\exp\left \{ -\frac{x ^{2}}{8\left (  u+\left ( \sigma \sqrt{n}  \right )^{-\alpha }x^{2-\alpha }  \right )   }  \right \}+\exp\left \{ -x^{2}  \right \} \nonumber \\
& \le & 3\exp\left \{ -\frac{x ^{2}}{8\left (  u+\left ( \sigma \sqrt{n}  \right )^{-\alpha }x^{2-\alpha }  \right )   }  \right \},\nonumber
\end{eqnarray}
which gives the desired inequality for $0\le x<\frac{\sigma \sqrt{n} }{2}.$\\
With arguments similar to that of \eqref{3.3} and \eqref{3.4}, we get for all $x>\frac{\sigma \sqrt{n} }{2}$,
\begin{eqnarray}
  I_{3}& = & \mathbb{P}\left ( \sum\limits _{i=1}^{n} \left ( X _{n_{0}+i }-\mu  \right )\ge \sigma \sqrt{n}\frac{x}{2} \right ) \nonumber \\
& \le  & 2\exp\left \{ -\frac{\left ( \frac{\sigma \sqrt{n}x }{2}  \right )^{2}  }{2\left ( u_{n}+\left ( \frac{\sigma \sqrt{n}x }{2}  \right )^{2-\alpha }   \right ) }  \right \} \nonumber \\
& = & 2\exp \left \{ -\frac{x^{2}  }{8\left ( \frac{u_{n}}{\sigma ^{2}n } +\left ( \sigma \sqrt{n}  \right ) ^{-\alpha }\left ( \frac{x}{2}  \right ) ^{2-\alpha  }    \right ) }  \right \} \nonumber \\
& \le & 2\exp \left \{ -\frac{  x^{2}  }{8\left ( u+\left ( \sigma \sqrt{n}  \right ) ^{-\alpha }x^{2-\alpha  }    \right )  }  \right \}\label{3.7}
\end{eqnarray}
and
\begin{eqnarray}
I_{4}\le\exp\left \{- \frac{x\sigma \sqrt{n} }{2}  \right \}.\label{3.8}
\end{eqnarray}
Combining \eqref{3.5}, \eqref{3.7} and \eqref{3.8}, we obtain for all $x>\frac{\sigma \sqrt{n} }{2}$,
\begin{eqnarray}
\mathbb{P}\left ( Z_{n_{0},n } \ge x \right )& \le & 2\exp\left \{ -\frac{x ^{2}}{8\left (  u+\left ( \sigma \sqrt{n}  \right )^{-\alpha }x^{2-\alpha }  \right )   }  \right \}+\exp\left \{ -\frac{x\sigma \sqrt{n} }{2}   \right \} \nonumber \\
& \le & 3\exp\left \{ -\frac{x ^{2}}{8\left (  u+\left ( \sigma \sqrt{n}  \right )^{-\alpha }x^{2-\alpha }  \right )   }  \right \}.\nonumber
\end{eqnarray}
This completes the proof of Theorem \ref{theorem 2.1}.  \hfill\qed
\subsection{Proof of Theorem \ref{theorem 2.2}}
From \eqref{3.1} and \eqref{3.2}, using Fuk-Nagaev's inequality  \cite{Nagaev}  for i.i.d.\ random variables,   we obtain for all $0< x\leq\frac{\sigma \sqrt{n} }{2},$
\begin{eqnarray}
I_{1}& \le & \exp\left \{ -\frac{ x^{2} \left ( 1-\frac{x }{\sigma \sqrt{n} }  \right )^{2}  }{\frac{1}{2}V^{2} }  \right \}+ \frac{ 2^{-(p+1)} n C_{p} }{\left (   \sqrt{n} \left ( x-\frac{x^{2} }{\sigma \sqrt{n} }  \right ) \right ) ^{p }}\nonumber\\
& \le  &\exp\bigg\{- \frac{  x^2 }{ 2 V^2  } \bigg\} +  \frac{ C_p}{2n^{(p-2)/2}x^p }, \label{6.1}
\end{eqnarray}
where $V^{2}$ and $C_{p }$ are defined in the theorem and the definition of $I_{i}(i=1,2,3,4) $ are shown in the proof of Theorem \ref{theorem 2.0}. By \eqref{3.4}, we have for all $0< x\leq\frac{\sigma \sqrt{n} }{2} $,
\begin{eqnarray}
I_2 \leq  \exp\Big\{   - \frac12 x^2 \Big\}  \leq  \frac{ C_p}{2n^{(p-2)/2}x^p }.\label{6.4}
\end{eqnarray}
Combining \eqref{6.1} and \eqref{6.4}, we have for all  $0<x<\frac{\sigma \sqrt{n} }{2},$
$$\mathbb{P}\left ( Z_{n_{0} ,n}\ge x  \right ) \leq  \exp\bigg\{- \frac{  x^2 }{ 2 V^2  } \bigg\} +  \frac{ C_p}{n^{(p-2)/2}x^p }, $$
which gives the desired inequality. When $x>\frac{\sigma \sqrt{n} }{2}$, it holds
$$\mathbb{P}\left ( Z_{n_{0},n }\ge x  \right )\le I_{3}+I_{4} .$$
Again by Fuk-Nagaev's inequality  \cite{Nagaev}  for i.i.d.\ random variables, we obtain for  all  $x>\frac{\sigma \sqrt{n} }{2}$,
\begin{eqnarray}
I_{3}& \le & \exp\left \{ -\frac{\left ( \frac{x}{2}  \right ) ^{2}  }{\frac{1}{2}V^{2} }  \right \}+\frac{2^{-(p+1)}nC_{p } }{\left ( \sqrt{n}  \frac{x}{2}  \right ) ^{p }   }\nonumber \\
& \le & \exp\bigg\{- \frac{  x^2 }{ 2 V^2  } \bigg\} +  \frac{ C_p}{2n^{(p-2)/2}x^p }.\label{6.2}
\end{eqnarray}
By an argument similar to that of  \eqref{3.8}, we have for all $x>\frac{\sigma \sqrt{n} }{2}$,
\begin{eqnarray}
I_{4}\le \exp\left \{- \frac{x\sigma \sqrt{n} }{2}  \right \}\le \frac{ C_p}{2n^{(p-2)/2}x^p }\label{6.3}
\end{eqnarray}
Combining \eqref{6.2} and \eqref{6.3}, we get for all $x>\frac{\sigma \sqrt{n} }{2}$,
\begin{eqnarray}
\mathbb{P}\left ( Z_{n_{0},n }\ge x  \right )& \le & \exp\bigg\{- \frac{  x^2 }{ 2 V^2  } \bigg\} +  \frac{ C_p}{n^{(p-2)/2}x^p } . \nonumber
\end{eqnarray}
This completes the proof of Theorem \ref{theorem 2.2}.  \hfill\qed
\subsection{Proof of Theorem \ref{theorem 2.5}}
Recall that $\ln Z_n = S_n + \ln W_n.$ It is easy to see that for all $x >0,$
\begin{eqnarray}
\mathbb{P} \left (\frac{1}{n} \ln \frac{Z_{n+n_{0}} }{Z_{n_{0} } }-\mu  \ge x\right )&\le& \mathbb{P}\left (\frac{\sum\limits_{i=1}^{n}\left ( X_{n_{0}+i }-\mu   \right ) }{n} \ge \frac{x}{2}    \right )+\mathbb{P}\left ( \frac{ \ln W_{n_{0},n }   }{n}\ge \frac{x}{2}   \right )\nonumber \\
&=:& K_{1}+K_{2} ,\label{9.4}
\end{eqnarray}
Using Markov's inequality and von Bahr-Esseen's inequality \cite{von Bahr} for i.i.d. random variables, we obtain for all $x> 0,$
\begin{eqnarray}
K_1 &\leq&\frac{ 2^p  }{ x^p}  \mathbb{E} \left | \sum_{i=1}^n\frac{X_{n_{0}+i } -\mu }{n}   \right | ^p \leq \frac{ 2^{p+1}  }{ x^p} \sum_{i=1}^n \mathbb{E} \left |\frac{X_{n_{0}+i } -\mu }{n}  \right | ^p\nonumber  \\
&= &\frac{ 2^{p+1} n\mathbb{E}|X -\mu  |^p }{ x^pn^{p} } =   2^{p+1} \mathbb{E}|X-\mu  |^p \frac{1}{ x^p n^{p-1} }.\label{2.1}
\end{eqnarray}
By an argument similar to that of  \eqref{3.8}, we have for all $x >  0,$
\begin{eqnarray}
K_2 \leq  \exp\Big\{   - \frac12 nx \Big\} \leq (2p)^pe^{-p} \frac{1}{  x^p n^{p-1} }.\label{2.2}
\end{eqnarray}
Combining \eqref{2.1} and \eqref{2.2} together, we obtain the desired inequality.  \hfill\qed
\subsection{Proof of Theorem \ref{theorem 2.6}}
 We first give a proof for \eqref{10}. From \eqref{3.1} and \eqref{3.2}, using Hoeffding's inequality \cite{Hoeffding} for i.i.d. random variables and the relation among the bounds of Hoeffding and Bernstein \cite{X.Fan1}, we get for all $0\le x<\frac{\sigma \sqrt{n} }{2},$
\begin{eqnarray}
I_{1}& \le &\exp \left \{ -\frac{x\left (  1-\frac{x}{\sigma \sqrt{n} } \right ) }{H}\left [ \left ( 1+\frac{\sigma ^{2} }{H\frac{\sigma }{\sqrt{n} }x\left (  1-\frac{x}{\sigma \sqrt{n} }\right )  }  \right ) \ln \left ( 1+\frac{H\frac{\sigma }{\sqrt{n} }x\left (  1-\frac{x}{\sigma \sqrt{n} }\right ) }{\sigma ^{2}}  \right )-1  \right ]   \right \}\nonumber \\
&\le& \exp \left \{ -\frac{\frac{\sigma ^{2}n }{H^{2} }   x^{2}\left ( 1-\frac{x}{\sigma \sqrt{n} }  \right )^{2}   }{2\left ( \frac{\sigma ^{2}n }{H^{2} }  +\frac{1}{3}\frac{\sigma \sqrt{n} }{H}  x\left ( 1-\frac{x}{\sigma \sqrt{n} }  \right ) \right )  }  \right \},\nonumber
\end{eqnarray}
where $I_{i}(i=1,2,3,4)$ are as defined in the proof of Theorem \ref{theorem 2.0}.
After some calculations, we get
\begin{eqnarray}
I_{1}& \le &\exp \left \{ -\frac{x}{2H } \left [ \left ( 1+\frac{2\sigma \sqrt{n}}{Hx }  \right )\ln \left ( 1+\frac{Hx }{2\sigma \sqrt{n}}  \right )  -1 \right ] \right \}\nonumber \\
&\le& \exp \left \{ -\frac{x^{2} }{8\left ( 1+\frac{Hx}{6\sigma \sqrt{n} }  \right ) }  \right \}.\label{7.1}
\end{eqnarray}
By Markov's inequality, the fact that $\mathbb{E}W_{n}=1 $ and the following inequality
\begin{eqnarray}
\frac{x}{2H } \left [ \left ( 1+\frac{2\sigma \sqrt{n}}{Hx }  \right )\ln \left ( 1+\frac{Hx }{2\sigma \sqrt{n}}  \right )  -1 \right ]\le \frac{x}{2H}\cdot \frac{Hx}{2\sigma \sqrt{n} }=\frac{x^{2} }{4\sigma \sqrt{n} },\nonumber
\end{eqnarray}
 we have for all $0\le x<\frac{\sigma \sqrt{n} }{2},$
\begin{eqnarray}
I_{2}&\le& \exp\left \{ -x^{2}  \right \}\le \exp \left \{ -\frac{x^{2} }{4\sigma \sqrt{n} } \right \}\nonumber  \\
&\le & \exp \left \{ -\frac{x}{2H } \left [ \left ( 1+\frac{2\sigma \sqrt{n}}{Hx }  \right )\ln \left ( 1+\frac{Hx }{2\sigma \sqrt{n}}  \right )  -1 \right ] \right \}\nonumber  \\
&\le & \exp \left \{ -\frac{x^{2} }{8\left ( 1+\frac{Hx}{6\sigma \sqrt{n} }  \right ) }  \right \}. \label{7.2}
\end{eqnarray}
Combining \eqref{7.1} and \eqref{7.2}, we obtain the desired inequality for all $0\le x<\frac{\sigma \sqrt{n} }{2}$.

Next, we give a proof for \eqref{11}. When $x>\frac{\sigma \sqrt{n} }{2},$ it holds
$$\mathbb{P}\left ( Z_{n_{0},n }\ge x  \right )\le I_{3}+I_{4} .$$
Again by Hoeffding's inequality \cite{Hoeffding} for i.i.d.  random variables and the relation among the bounds of Hoeffding and Bernstein \cite{X.Fan1}, we obtain for all $x>\frac{\sigma \sqrt{n} }{2},$
\begin{eqnarray}
I_{3}& \le &\exp \left \{ -\frac{x}{2H } \left [ \left ( 1+\frac{2\sigma \sqrt{n}}{Hx }  \right )\ln \left ( 1+\frac{Hx }{2\sigma \sqrt{n}}  \right )  -1 \right ] \right \}\nonumber \\
&\le & \exp \Bigg \{ -\frac{x^{2} }{8\big( 1+\frac{Hx}{6\sigma \sqrt{n} }  \big) }  \Bigg\}.\label{2.5}
\end{eqnarray}
By an argument similar to that of  \eqref{2.8}, we have for all $x>\frac{\sigma \sqrt{n} }{2},$
\begin{eqnarray}
I_{4}\le \exp\left \{- \frac{x\sigma \sqrt{n} }{2}  \right \}.\label{2.6}
\end{eqnarray}
Combining \eqref{2.5} and \eqref{2.6}, we get for all  $x>\frac{\sigma \sqrt{n} }{2},$
\begin{eqnarray}
\mathbb{P}\left ( Z_{n_{0},n }\ge x  \right ) &\le& \exp \left \{ -\frac{x}{2H } \left [ \left ( 1+\frac{2\sigma \sqrt{n}}{Hx }  \right )\ln \left ( 1+\frac{Hx }{2\sigma \sqrt{n}}  \right )  -1 \right ] \right \}+\exp\left \{- \frac{x\sigma \sqrt{n} }{2}  \right \}\nonumber \\
&\le&\exp \Bigg \{ -\frac{x^{2} }{8\left ( 1+\frac{Hx}{6\sigma \sqrt{n} }  \right ) }  \Bigg \}+\exp\Bigg \{- \frac{x\sigma \sqrt{n} }{2}  \Bigg \},\nonumber
\end{eqnarray}
which gives the desired inequality.\hfill\qed
\subsection{Proofs of Theorem \ref{theorem 2.8}}
 From \eqref{9.4}, using Rio's inequality \cite{Rio} for i.i.d.\ random variables, we get for all $x\in [0,2\left ( H_{2}-H_{1}   \right ))  $,
\begin{eqnarray}
K_{1}&=&\mathbb{P}\left ( \frac{\sum_{i=1}^{n}\left ( X_{n_{0}+i }-\mu   \right )}{H_{2}-H_{1}  } \ge n\cdot \frac{x}{2\left ( H_{2}-H_{1} \right ) }    \right )\nonumber \\
&\le &\exp \Big \{ -n\max \left ( \psi _{1}\left ( x \right ), \psi _{2}\left ( x \right )   \right ) \Big \}\nonumber \\
&\le &\left ( 1-\frac{x}{2\left ( H_{2}-H_{1} \right )}  \right )^{\frac{nx}{H_{2}-H_{1}} \left ( 1-\frac{x }{4\left ( H_{2}-H_{1} \right ) }  \right ) },\label{0.6}
\end{eqnarray}
where $\psi _{1}\left ( x \right ) $ and $\psi _{2}\left ( x \right )$ are defined as \eqref{0.8} and $K_{i}(i=1,2)$ are as defined in the proof of Theorem \ref{theorem 2.5}. Next, by Markov's inequality and the fact that $\mathbb{E}W_{n}=1 $, we have for all $x\ge 0,$
\begin{eqnarray}
K_{2}\le \exp\left \{ -\frac{1}{2} nx \right \}. \label{0.7}
\end{eqnarray}
Combining \eqref{0.6} and \eqref{0.7}, we obtain for all $x\in [0,2\left ( H_{2}-H_{1}   \right ))$,
\begin{eqnarray}
\mathbb{P} \left (\frac{1}{n} \ln \frac{Z_{n+n_{0}} }{Z_{n_{0} } }-\mu  \ge x\right )&\le&\exp \Big \{ -n\max \left ( \psi _{1}\left ( x \right ), \psi _{2}\left ( x \right )   \right ) \Big \} +\exp \left \{ -\frac{1}{2}nx  \right \}\nonumber \\
&\le&\left ( 1-\frac{x}{2\left ( H_{2}-H_{1} \right )}  \right )^{\frac{nx}{H_{2}-H_{1}} \left ( 1-\frac{x }{4\left ( H_{2}-H_{1} \right ) }  \right ) } +\exp \left \{ -\frac{1}{2}nx  \right \}.\nonumber
\end{eqnarray}
This completes the proof of Theorem \ref{theorem 2.8}.  \hfill\qed

\section*{Acknowledgements}
Xu would like to thank Xiequan Fan for his helpful suggestions.


\begin{thebibliography}{00}{\footnotesize


\bibitem{Afanasyev} V.I. Afanasyev, C. B{\"o}inghoff, G. Kersting, V.A. Vatutin, Limit theorems for weakly subcritical branching processes in random environment, J. Theoret. Probab. 25 (3) (2012) 703-732.
\bibitem{Boinghoff} V.I. Afanasyev, C. B{\"o}inghoff, G. Kersting, V.A. Vatutin, Conditional limit theorems for intermediately subcritical branching processes in random environment, Ann. Inst. Henri Poincar{\'e} Probab. Stat. 50 (2) (2014) 602-627.
\bibitem{Bansaye} V. Bansaye, J. Berestycki, Large deviations for branching processes in random environment, Markov Process. Related Fields 15 (4) (2009) 493-524.
\bibitem{Bansaye2} V. Bansaye, C. B{\"o}inghoff, Upper large deviations for branching processes in random environment with heavy tails, Electron. J. Probab. 16 (69) (2011) 1900-1933.
\bibitem{Bor} A.A. Borovkov, Estimates for the distribution of sums and maxima of sums of random variables when the Cram\'er  condition is not satisfied, Sib. Math. J. 41 (2000) 811-848.
\bibitem{Boinghoff2} C. B{\"o}inghoff, G. Kersting, Upper large deviations of branching processes in a random environemnt - offspring distributions with geomertrically bounded tails, Stochastic Process. Appl. 120 (10) (2010) 2064-2077.
\bibitem{Boinghoff3} C. B{\"o}inghoff, Limit theorems for strongly and intermediately supercritical branching processes in random environment with linear fractional offspring distributions, Stochastic Process. Appl. 124 (11) (2014) 3553-3577.
\bibitem{Bernstein} S.N. Bernstein, The Theory of Probabilities, Moscow: Leningrad, 1946.
\bibitem{DDF19} J. Dedecker, P. Doukhan,  X. Fan,  Deviation inequalities for separately Lipschitz functionals of composition of random functions.
J.  Math. Anal.  Appl. 479(2) (2019) 1549--1568.
\bibitem{D99} V.H. De la Pe\~{n}a,  A general class of exponential inequalities for martingales and ratios, Ann. Probab. 27 (1999) 537-564.
\bibitem{FGL12} X. Fan, I. Grama, Q. Liu, Cram\'{e}r large deviation expansions for martingales under Bernstein's condition, Stochastic Process. Appl. 123 (2013) 3919-3942.
\bibitem{X.Fan} X. Fan, I. Grama, Q. Liu, Deviation inequalities for martingales with applications, J. Math. Anal. Appl. 448 (2017) 538-566.
\bibitem{X.Fan1} X. Fan, I. Grama, Q. Liu, Hoeffding's inequality for supermartingales, Stochastic Process. Appl. 122 (10) (2012) 3545-3559.
\bibitem{Fan} X. Fan, H. Hu, Q. Liu, Uniform Cram\'er moderate deviations and Berry-Esseen bounds for a supercritical branching process in a random environment. Front. Math. China. 15 (5) (2020) 891-914.
\bibitem{Grama} I. Grama, Q. Liu, M. Miqucu, Berry-Esseen's bound and Cram{\'e}r's large deviations for a supercritical branching process in a random environment. Stochastic Process. Appl. 127 (2017) 1255-1281.
\bibitem{Huang} C. Huang, Q. Liu, Moments, moderate and large deviations for a branching process in a random environment, Stochastic Process. Appl. 122 (2) (2012) 522-545.
\bibitem{Hoeffding} W. Hoeffding, Probalility inequalities for sums of bounded random variables, J. Amer. Statist. Assoc. 58 (1963) 13-30.
\bibitem{Kozlov} M.V. Kozlov, On large deviations of branching processes in a random environment: geometric distribution of descendants, Discrete Math. Appl. 16 (2) (2006) 155-174.

\bibitem{Nakashima} M. Nakashima, Lower deviations of branching processes in random environment with geometrical offspring distributions, Stochastic Process. Appl. 123 (9) (2013) 3560-3587.
\bibitem{Nagaev} S.V. Nagaev, Large deviations of sums of independent random variables, Ann. Probab. 7 (1979) 745-789.
\bibitem{NV} S.V. Nagaev, V. Vakhtel, Probability inequalities for a critical Galton--Watson process, Theory Probab. Appl. 50 (2) (2006) 225C247.
\bibitem{Rio} E. Rio, On McDiarmid's concentration inequality, Electron. Commun. Probab. 18 (2013) 1-11.
\bibitem{Smith} W.L. Smith, W.E. Wilkinson, On branching processes in random environments, Ann. Math. Stat. 40 (3) (1969) 814-827.
\bibitem{Vatutin} V.A. Vatutin, A refinement of limit theorems for the critical branching processes in random environment, in: Workshop on Branching Processes and their Applications, in: Lect. Notes Stat. Proc., vol. 197, Springer,Berlin (2010) 3-19.
\bibitem{von Bahr} B. von Bahr, C.G. Esseen, Inequlities for the $r$th absolute moment of a sum of random variables, $1 \leqq r \leqq 2$, Ann. Math. Stat. 36 (1) (1965) 299-303.
\bibitem{WL17}  Y. Wang,  Q. Liu,  Limit theorems for a supercritical branching process
with immigration in a random environment, Sci. China Math. 60(12) (2017) 2481-2502.
}
\end{thebibliography}
\end{document}